\documentclass[a4, 12pt]{amsart}
\usepackage{amssymb}
\usepackage{amstext}
\usepackage{amsmath}
\usepackage{amscd}
\usepackage{latexsym}
\usepackage{amsfonts}

\theoremstyle{plain}
\newtheorem{thm}{Theorem}[section]
\newtheorem*{thm*}{Theorem}
\newtheorem*{cor*}{Corollary}

\newtheorem{prop}[thm]{Proposition}
\newtheorem{lem}[thm]{Lemma}
\newtheorem{cor}[thm]{Corollary}

\newtheorem*{claim*}{Claim}

\theoremstyle{definition}

\newtheorem{ex}[thm]{Example}

\theoremstyle{remark}

\numberwithin{equation}{thm}

\def\rank{\rm{rank}}
\def\a{\mathfrak a}
\def\b{\mathfrak b}
\def\c{\mathfrak c}
\def\e{\mathrm e}
\def\m{\mathfrak m}

\def\p{\mathfrak p}
\def\q{\mathfrak q}

\def\Z{\Bbb Z}

\newcommand{\rmG}{\rm{G}}

\newcommand{\rmR}{\rm{R}}

\newcommand{\fkm}{\mathfrak{m}}

\newcommand{\fkp}{\mathfrak{p}}

\def\Ass{\rm{Ass}}

\begin{document}

\setlength{\baselineskip}{20pt}

\title{The structure of Sally modules of rank one}
\author{Shiro Goto}
\address{Department of Mathematics, School of Science and Technology, Meiji University, 1-1-1 Higashi-mita, Tama-ku, Kawasaki 214-8571, Japan}
\email{goto@math.meiji.ac.jp}
\author{Koji Nishida}
\address{Department of Mathematics and Informatics, Graduate School of Science and Technology, Chiba University, 1-33 Yatoi-cho, Inage-ku, Chiba-shi, 263 Japan}
\email{nishida@math.s.chiba-u.ac.jp}
\author{Kazuho Ozeki}
\address{Department of Mathematics, School of Science and Technology, Meiji University, 1-1-1 Higashi-mita, Tama-ku, Kawasaki 214-8571, Japan}
\email{kozeki@math.meiji.ac.jp}
\thanks{
\endgraf
}
\maketitle
\begin{abstract}
A complete structure theorem of Sally modules of $\fkm$-primary ideals $I$ in a Cohen-Macaulay local ring 
$(A, \m)$ satisfying the equality $\e_1(I)=\e_0(I)-\ell_A(A/I)+1$  is given, 
where $\e_0(I)$ and $\e_1(I)$ denote the first two Hilbert coefficients of $I$.
\end{abstract}
\section{Introduction}
This paper aims to give a structure theorem of Sally modules of rank one.

Let $A$ be a Cohen-Macaulay local ring with the maximal ideal $\m$ and $d=\dim A >0$. We assume the residue class field $k=A/\fkm$ of $A$ is infinite. Let $I$ be an $\fkm$-primary ideal in $A$ and choose a minimal reduction $Q=(a_1, a_2, \cdots, a_d)$ of $I$. Then we have integers $\{e_i=\e_i(I)\}_{0 \leq i \leq d}$ such that the equality
$$\ell_A(A/I^{n+1})={e}_0\binom{n+d}{d}-{e_1}\binom{n+d-1}{d-1}+\cdots+(-1)^d{e}_d$$
holds true for all $n \gg 0$. Let $$R = \mathrm{R}(I) := A[It]~~ \ \operatorname{and}~~ \ T= \mathrm{R}(Q):= A[Qt]~~\subseteq~~A[t]$$ denote, respectively, the Rees algebras of $I$ and $Q$, where $t$ stands for an indeterminate over $A$. We put $$R' = \mathrm{R}'(I) := A[It, t^{-1}] \ \ \operatorname{and}  \ \  G= \mathrm{G}(I) :=  R'/t^{-1}R'~~\cong~~\bigoplus_{n \geq 0}I^n/I^{n+1}.$$ 
Let $B=T/\m T$, which is the polynomial ring with $d$ indeterminates over the field $k$. Following W. V. Vasconcelos \cite{V}, we then define
 $$\mathrm{S}_Q(I)= IR/IT$$ and call it the Sally module of $I$ with respect to $Q$. We notice that the Sally module $S=\mathrm{S}_Q(I)$ is a finitely generated graded $T$-module, since $R$ is a module-finite extension of the graded ring $T$.

The Sally module $S$ was introduced by W. V. Vasconcelos \cite{V}, 
where he gave an elegant review, in terms of his $Sally$ module, of the works \cite{S1, S2, S3}
of J. Sally about the structure of $\m$-primary ideals $I$ with interaction to the structure
of the graded ring $G$ and the Hilbert coefficients $e_i$'s of $I$.

As is well-known, we have the inequality  (\cite{N}) $$e_1 \geq e_0 - \ell_A(A/I)$$ and  C. Huneke  \cite{H} showed that $e_1 = e_0 - \ell_A(A/I)$ if and only if $I^2 = QI$. When this is the case, both the graded rings $G$ and $\operatorname{F}(I) = \bigoplus_{n \geq 0}I^n/\fkm I^n$ are Cohen-Macaulay, and the Rees algebra $R$ of $I$  is also a Cohen-Macaulay ring, provided $d \geq 2$. Thus, the ideals $I$ with $e_1 = e_0 - \ell_A(A/I)$ enjoy very nice properties. The reader may consult with the recent work of Wang \cite{W}, which establishes the ubiquity of ideals $I $ with $I^2 =QI$.

J. Sally \cite{S3}  firstly investigated the second border, that is the ideals $I$ satisfying the equality 
$$e_1 = e_0 - \ell_A(A/I) + 1$$ and gave several  very important results. Among them, one can find the following characterization of ideals $I$ with 
$e_1 = e_0 - \ell_A(A/I) + 1$ and $e_2 \ne 0$, where $B(-1)$ stands for the graded $B$-module 
whose grading is given by $[B(-1)]_n = B_{n-1}$ for all $n \in \Z$. 
The reader may also consult with \cite{CPP} and \cite{P} for further ingenious use of Sally modules.

\begin{thm}[Sally \cite{S3} , Vasconcelos \cite{V}]\label{Sally}
The following three conditions are equivalent to each other.
\begin{itemize}
\item[$(1)$] $S \cong B(-1)$ as graded $T$-modules.
\item[$(2)$] $e_1=e_0-\ell_A(A/I)+1$ and if $d \geq 2$, $e_2 \ne 0$.
\item[$(3)$] $I^3 = QI^2$ and $\ell_A(I^2/QI) = 1$.
\end{itemize}
When this is the case, the following assertions hold true.
\begin{itemize}
\item[(i)] $e_2 = 1$, if $d \geq 2$.
\item[(ii)] $e_i = 0$ for all $3 \leq i \leq d$.
\item[(iii)] $\operatorname{depth} G \geq d-1$.
\end{itemize}
\end{thm}

This beautiful theorem says, however, nothing about the case where $e_2=0$. It seems natural to ask what happens, when $e_2 = 0$,  on the ideals $I$ which satisfy the equality $e_1 = e_0 - \ell_A(A/I) +1$. This long standing question has motivated the recent research \cite{GNO}, where the authors gave several partial answers to the question. The present research is a continuation of \cite{GNO, S3, V} and aims at a simultaneous understanding of 
the structure of Sally modules of ideals $I$ which satisfy the equality 
$e_1=e_0- \ell_A(A/I)+1$. 

Let us now state our own result. The main result of this paper is the following Theorem \ref{MainTheorem}, which contains Theorem \ref{Sally} of Sally--Vasconcelos as the case where $c = 1$. Our contribution in Theorem \ref{MainTheorem} is the implication $(1) \Rightarrow (3)$, the proof of which is based on the new result that the equality $I^3=QI^2$ holds true if $e_1=e_0 - \ell_A(A/I) +1$ (cf. Theorem \ref{I3=QI2}).

\begin{thm}\label{MainTheorem}
The following three conditions are equivalent to each other.
\begin{itemize}
\item[$(1)$] $e_1=e_0 - \ell_A(A/I) +1$. 
\item[$(2)$] $\fkm S = (0)$ and $\operatorname{rank}_BS = 1$.
\item[$(3)$] $S \cong (X_1,X_2, \cdots ,X_c)B$
as graded $T$-modules for some $0 < c \leq d$,
where $\{X_i\}_{1 \leq i \leq c}$ are linearly independent linear forms of the polynomial ring $B$.
\end{itemize}
When this is the case, $c = \ell_A(I^2/QI)$ and $I^3=QI^2$, and the following assertions hold true.
\begin{itemize}
\item[(i)] $\operatorname{depth} G \geq d-c$ and $\operatorname{depth}_TS=d-c+1$.
\item[(ii)] $\operatorname{depth} G = d-c$, if $c \geq 2$.
\item[(iii)] Suppose $c < d$. Then 
$$\ell_A(A/I^{n+1})=e_0 \binom{n+d}{d}-e_1 \binom{n+d-1}{d-1}+\binom{n+d-c-1}{d-c-1}$$ for all $n \geq 0$.
Hence
\[ e_i =  \left\{
\begin{array}{rl}
0 & \quad \mbox{if $i \neq c+1 $,} \\
(-1)^{c+1} & \quad \mbox{if $i =c+1$}
\end{array}
\right.\]
for $2 \leq i \leq d$.
\item[(iv)] Suppose $c = d$. Then 
$$\ell_A(A/I^{n+1})=e_0 \binom{n+d}{d}-e_1 \binom{n+d-1}{d-1}$$
for all $n \geq 1$.
Hence $e_i=0$  for $2 \leq i \leq d$.
\end{itemize}
\end{thm}

Thus Theorem \ref{MainTheorem} settles a long standing problem, although the structure of ideals $I$ with $e_1 = e_0 - \ell_A(A/I) + 2$ or the structure of Sally modules $S$ with  $ \fkm S = (0)$ and $\operatorname{rank}_BS=2$ remains unknown.

Let us now briefly explain how this paper is organized. 
We shall prove Theorem \ref{MainTheorem} in Section 3. 
In Section 2 we will pick up from the paper \cite{GNO} some auxiliary results on Sally modules, all of which  are  known, but let us  note them for the sake of the reader's convenience. 
In Section 4 we shall discuss two consequences of Theorem \ref{MainTheorem}. The results are more or less  known by \cite{GNO, S3, V}. However, thanks to Theorem \ref{MainTheorem}, not only the statements of the results but also the proofs are substantially simplified, so that we would like to note the improved statements, and would like to indicate a brief proof of Theorem \ref{Sally} as well. 
In Section 5 we will construct one example in order to see the ubiquity of ideals $I$ which satisfy condition (3) in Theorem \ref{MainTheorem}. We will show that, for given integers $0 < c \leq d$, there exists an $\fkm$-primary ideal $I$  in a certain Cohen-Macaulay local ring $(A, \fkm)$ such that  $$d = \operatorname{dim}A, \ \ e_1 = e_0 - \ell_A(A/I) + 1, \ \ \operatorname{and} \ \ c = \ell_A(I^2/QI)$$ for some reduction $Q = (a_1, a_2, \cdots, a_d) $ of $I$.

In what follows, unless otherwise specified, let  $(A, \m)$ be a Cohen-Macaulay local ring with $d=\dim A >0$. We assume that the field $k=A/ \m$ is infinite. Let $I$ be an $\m$-primary ideal in $A$ and let $S$ be the Sally module of $I$ with respect to a minimal reduction $Q = (a_1, a_2, \cdots, a_d)$ of $I$. We put $R= A[It], T=A[Qt], R' = A[It, t^{-1}]$, and $G= R'/t^{-1}R'$. 
Let $$\tilde{I} = \bigcup_{n \geq 1}[I^{n+1} :_AI^n]$$ denote the Ratliff-Rush closure of $I$, which is the largest $\fkm$-primary ideal in $A$ such that $I \subseteq \tilde{I}$ and  $\e_i(\tilde{I}) = e_i$ for all $0 \leq i \leq d$ (cf. \cite{RR}). We denote by $\mu_A(*)$ the number of generators.


\section{Auxiliary results}

In this section let us  firstly summarize some known results on Sally modules, which we need throughout this paper. See \cite {GNO} and \cite{V} for the detailed proofs.

The first two results are basic facts on Sally modules developed by Vasconcelos \cite{V}.

\begin{lem}\label{fact1}
The following assertions hold true.
\begin{itemize}
\item[$(1)$] $\m^{\ell} S = (0)$ for  integers $\ell \gg 0$.
\item[$(2)$]The homogeneous components $\{ S_n \}_{n \in \Z}$ of the graded $T$-module $S$ are given by
\[ S_n \cong  \left\{
\begin{array}{rl}
(0) & \quad \mbox{if $n \leq 0 $,} \\
I^{n+1}/IQ^n & \quad \mbox{if $n \geq 1$.}
\end{array}
\right.\]
\item[$(3)$] $S=(0)$ if and only if $I^2=QI$.
\item[$(4)$] Suppose that $S \ne (0)$ and put $V = S/MS$, where $M = \fkm T + T_+$ is the graded maximal ideal in $T$. Let $V_n$~$(n \in \Z)$  denote the homogeneous component of the finite-dimensional graded $T/M$-space $V$ with degree $n$ and put $\Lambda = \{n \in \Z \mid V_n \ne (0) \} $. Let $q = \max \Lambda$. Then we have $\Lambda =\{1, 2, \cdots, q\}$ and $\mathrm{r}_Q(I) = q+1$, where $\mathrm{r}_Q(I)$ stands for  the reduction number of $I$ with respect to $Q$.   
\item[$(5)$] $S = TS_1$ if and only if $I^3 = QI^2$. 
\end{itemize}
\end{lem}

\begin{proof}
See \cite[Lemma 2.1]{GNO}.
\end{proof}

\begin{prop}\label{fact2}
Let $\fkp = \m T$. Then the following assertions hold true.
\begin{itemize}
\item[$(1)$] ${\Ass}_TS \subseteq \{ \fkp \}$. Hence $\dim_TS=d$, if $S\ne (0)$. 
\item[$(2)$] $\ell_A(A/I^{n+1})={e}_0\binom{n+d}{d}-({e}_0-\ell_A(A/I)){\cdot}\binom{n+d-1}{d-1}-\ell_A(S_n)$ for all $n \geq 0$.
\item[$(3)$] We have $e_1 =e_0 - \ell_A(A/I) + \ell_{T_\fkp}(S_\fkp)$. Hence
$e_1 = e_0 - \ell_A(A/I) + 1$ if and only if $\m S = (0)$ and $\operatorname{rank}_BS = 1$. 
\item[$(4)$] Suppose that $S\ne (0)$. Let $s = \operatorname{depth}_TS$. Then  $\operatorname{depth} G = s-1$ if $s < d$. $S$ is a Cohen-Macaulay $T$-module if and only if  $\operatorname{depth} G \geq d-1$.
\end{itemize}
\end{prop}

\begin{proof}
See \cite[Proposition 2.2]{GNO}.
\end{proof}

Combining Lemma \ref{fact1} $(3)$ and Proposition \ref{fact2}, 
we readily get the following results of Northcott \cite{N} and Huneke \cite{H}.

\begin{cor}[\cite{H, N}]\label{huneke}
We have $e_1 \geq e_0 - \ell_A(A/I)$. 
The equality $e_1 = e_0 - \ell_A(A/I)$ holds true if and only if $I^2 = QI$.
When this is the case, $e_i = 0$ for all $2 \leq i \leq d$.
\end{cor}

The following result is one of the keys for our proof of Theorem \ref{MainTheorem}.

\begin{thm}[\cite{GNO}]\label{rank1}
The following conditions are equivalent.
\begin{itemize}
\item[$(1)$] $\fkm S=(0)$ and ${\rank}_BS=1$.
\item[$(2)$] $S \cong \a$ as graded $T$-modules for some graded ideal $\a~(\ne B)$ of $B$.
\end{itemize}
\end{thm}

\begin{proof}
See \cite[Theorem 2.4]{GNO}.
\end{proof}

The following result is also due to \cite{GNO}, which will enable us to reduce the proof of Theorem \ref{MainTheorem} to the proof of the fact that $I^3 = QI^2$ if $e_1 = e_0 - \ell_A(A/I) + 1$.

\begin{prop}[\cite{GNO}]\label{Theorem2}
Suppose $e_1 = e_0 - \ell_A(A/I) + 1$ and $I^3 =QI^2$. Let $c = \ell_A(I^2/QI)$. Then the following assertions hold true.
\begin{enumerate}
\item[$(1)$] $0 < c \leq d$ and $\mu_B(S) = c$.
\item[$(2)$] $\operatorname{depth} G \geq d-c$ and $\operatorname{depth}_BS = d-c + 1$. 
\item[$(3)$] $\operatorname{depth} G =d-c$, if $c \geq 2$.
\item[$(4)$] Suppose $c < d$. Then 
$\ell_A(A/I^{n+1})=e_0\binom{n+d}{d} -e_1\binom{n+d-1}{d-1}+\binom{n+d-c-1}{d-c-1}$ for all $n \geq 0$. Hence 
\[e_i =\left\{
\begin{array}{
rl}
0 \ \ \ \ \ & \quad \mbox{if $i \ne c+1$} \\
(-1)^{c+1} & \quad \mbox{if $i = c+1$}
\end{array}\right.\]
for $2 \leq i \leq d$.
\item[$(5)$] Suppose $c=d$. Then 
$\ell_A(A/I^{n+1})=e_0\binom{n+d}{d} -e_1\binom{n+d-1}{d-1}$ for all $n \geq 1$. Hence
$e_i = 0$ for $2 \leq i \leq d$. 
\end{enumerate}
\end{prop}

\begin{proof}
See \cite[Corollary 2.5]{GNO}.
\end{proof}

The following result might be known. However, since we can find no good references, let us include a brief proof.

\begin{prop}\label{2.6}
Let $Q \subseteq I \subseteq J$ be ideals in a commutative ring $A$. Assume that $J = I +(h)$ for some $h \in A$. Then $I^3 = QI^2$, if $J^2 = QJ$.
\end{prop}
\begin{proof}
Since $hI \subseteq J^2=Q J=QI+Qh$,
for each $i \in I$ there exist $j \in QI$ and $q \in Q$ such that $hi=j+qh$.
Hence $h(i-q) =j \in QI$.
On the other hand, we have
$(i-q)I^2 \subseteq (i-q) J^2=(i-q)(QI+Qh)=(i-q)QI+jQ \subseteq QI^2$,
because $i-q \in I$ and $j \in QI$.
Thus $(i-q)I^2 \subseteq QI^2$, so that we have $iI^2 \subseteq QI^2$ for all $i \in I$.
Hence $I^3=QI^2$.
\end{proof}


\section{Proof of Theorem \ref{MainTheorem}}
The purpose of this section is to prove Theorem \ref{MainTheorem}. See Proposition \ref{fact2} (3) for the equivalence of conditions (1) and (2) in Theorem \ref{MainTheorem}. The implication $(3) \Rightarrow (2)$ is clear. So, we must show the implication $(1) \Rightarrow (3)$ together with the last assertions in Theorem \ref{MainTheorem}. Suppose that $e_1 = e_0 -\ell_A(A/I) + 1$. Then, thanks to Theorem \ref{rank1}, we get an isomorphism $$\varphi : S \to \a$$ of graded $B$-modules,
where $\a \subsetneq B$ is a graded ideal of $B$. Notice that once we are able to  show $I^3 = QI^2$, the last assertions of Theorem \ref{MainTheorem} readily follow from Proposition \ref{Theorem2}. On the other hand, since $\a \cong S = BS_1$ (cf. Lemma \ref{fact1} (5)), the ideal  $\a$ of $B$ is generated by linearly independent linear forms $\{X_i\}_{1 \leq i \leq c}$ $(0 < c \leq d)$ of $B$ and so, the implication $(1) \Rightarrow (3)$ in Theorem \ref{MainTheorem} follows. We have $c = \ell_A(I^2/QI)$, because $\a_1 \cong S_1 = I^2/QI$ (cf. Lemma \ref{fact1} (2)). Thus  our Theorem \ref{MainTheorem} has been proven modulo the following theorem, which follows also,  in the case where $d \leq 2$, from a result of  M. Rossi \cite[Corollary 1.5]{R}.


\begin{thm}\label{I3=QI2}
Suppose that $e_1 = e_0 -\ell_A(A/I) + 1$. Then $I^3 = QI^2$.
\end{thm}

\begin{proof}
We proceed by induction on $d$. Suppose that $d = 1$. Then $S$ is $B$-free of rank one (recall that the $B$-module $S$ is torsionfree; cf. Proposition \ref{fact2} $(1)$) and so, since $S_1 \ne (0)$ (cf. Lemma \ref{fact1} (3)), $S \cong B(-1)$ as graded $B$-modules. Thus $I^3=QI^2$ by Lemma \ref{fact1} (5).

Let us assume that $d \geq 2$ and that our assertion holds true for $d-1$. Since the field $k=A/\fkm$ is infinite, without loss of generality we may assume that $a_1$ is a superficial element of $I$. Let $$\overline{A}=A/(a_1), \ \ \overline{I}=I/(a_1), \ \operatorname{and} \ \ \overline{Q}=Q/(a_1).$$ We then have $\e_i(\overline{I})=e_i$ for all $ 0 \leq i \leq d-1$, whence $$\e_1(\overline{I}) = \e_0(\overline{I}) - \ell_{\overline{A}}(\overline{A}/\overline{I}) + 1.$$ Therefore the hypothesis of induction on $d$  yields $\overline{I}^3 = \overline{Q}\,\overline{I}^2$. Hence, because the element $a_1t$ is a nonzerodivisor on $G$ if $\operatorname{depth} G > 0$, we have $I^3 = QI^2$ in that case.

Assume that $\operatorname{depth}G = 0$. Then, thanks to Sally's technique (\cite{S3}), we also have $\operatorname{depth} {\rmG}(\overline{I}) = 0$. Hence $\ell_{\overline{A}}(\overline{I}^2/\overline{Q}~\overline{I}) = d - 1$ by Proposition \ref{Theorem2} (2), because $\e_1(\overline{I}) = \e_0(\overline{I}) - \ell_{\overline{A}}(\overline{A}/\overline{I}) + 1.$ Consequently, $\ell_A(S_1) = \ell_A(I^2/QI) \geq d-1$, because $\overline{I}^2/\overline{Q}~\overline{I}$ is a homomorphic image of $I^2/QI$. Let us take an isomorphism $$\varphi : S \to \a$$ of graded $B$-modules,
where $\a \subsetneq B$ is a graded ideal of $B$. Then, since  $$\ell_A(\a_1) = \ell_A(S_1) \geq d-1,$$ the ideal $\a$ contains $d-1$ linearly independent linear forms, say $X_1, X_2, \cdots , X_{d-1}$ of $B$, which we enlarge to a basis $X_1, \cdots , X_{d-1}, X_d$ of $B_1$. Hence $$B = k[X_1, X_2, \cdots, X_{d}],$$ so that the ideal $\a/(X_1, X_2, \cdots, X_{d-1})B$ in the polynomial ring  $$B/(X_1, X_2, \cdots, X_{d-1})B = k[X_d]$$ is principal. If $\a = (X_1, X_2, \cdots, X_{d-1})B$, then $I^3 = QI^2$ by Lemma \ref{fact1} (5), since $S = BS_1$. However, because $\ell_A(I^2/QI) = \ell_A(\a_1) = d-1$, we have $\operatorname{depth} G \geq 1$ by Proposition \ref{Theorem2} (2), which is impossible. Therefore $\a/(X_1, X_2, \cdots, X_{d-1})B \ne (0)$, so that we have $$\a = (X_1, X_2, \cdots, X_{d-1}, X_d^\alpha)B$$ for some $\alpha \geq 1$. Notice that $\alpha = 1$ or $\alpha = 2$ by Lemma \ref{fact1} (4). 
We must show that $\alpha=1$.

Assume that $\alpha = 2$. Let us write, for each $1 \leq i \leq d$, $X_i = \overline{b_it}$ with $b_i \in Q$, where $\overline{b_it}$ denotes the image of $b_it \in T$ in $B=T/\fkm T$. Then $\a = (\overline{b_1t}, \overline{b_2t}, \cdots, \overline{b_{d-1}t}, \overline{(b_dt)^2})$. Notice that  
$$Q = (b_1, b_2, \cdots, b_d),$$
because $\{X_i\}_{1 \leq i \leq d}$ is a $k$-basis of $B_1$.
We now choose elements $f_i \in S_1$ for $1 \leq i \leq d-1$ and $f_d \in S_2$ so that $\varphi(f_i)=X_i$ for $1 \leq i \leq d-1$ and $\varphi(f_d) = X_d^2$. Let $z_i \in I^2$ for $1 \leq i \leq d-1$ and $z_d \in I^{3}$ such that
$\{f_i\}_{1 \leq i \leq d-1}$ and $f_d$ are, respectively, 
the images of $\{z_it\}_{1 \leq i \leq d-1}$ and $z_dt^2$ in $S$. We now consider the relations $X_if_1=X_1f_i$ in $S$ for $1 \leq i \leq d-1$ and $X_d^2f_1 = X_1f_d$, that is
$$b_i z_1-b_1 z_i \in Q^2I$$ for $1 \leq i \leq d-1$ and
$$b_d^2 z_1-b_1 z_d \in Q^{3}I.$$
Notice that $$Q^{3}=b_1Q^2+(b_2,b_3,\cdots,b_{d-1})^2{\cdot}(b_2,b_3,\cdots,b_{d}) + b_d^2Q$$ and write $$b_d^2 z_1-b_1 z_d=b_1 \tau_1+\tau_2+b_d^2 \tau_3$$
with $\tau_1 \in Q^2I$, $\tau_2 \in (b_2,b_3,\cdots,b_{d-1})^2{\cdot}(b_2,b_3,\cdots,b_d)I$,
and $\tau_3 \in QI$.
Then $$b_d^2 (z_1-\tau_3)=b_1(\tau_1+z_d)+\tau_2 \in (b_1)+(b_2,b_3,\cdots,b_{d-1})^2.$$
Hence $z_1-\tau_3 \in (b_1)+(b_2,b_3,\cdots,b_{d-1})^2$, 
because the sequence $b_1,b_2,\cdots,b_d$ is $A$-regular.
Let $z_1-\tau_3=b_1h+h'$ with $h \in A$ and $h' \in (b_2,b_3,\cdots,b_{d-1})^2$. Then since  $$b_1[b_d^2h - (\tau_1 +z_d)] = \tau_2 - b_d^2h' \in (b_2, b_3, \cdots, b_d)^3,$$ we have $b_d^2h -(\tau_1 + z_d) \in (b_2, b_3, \cdots, b_d)^3$, whence  $b_d^2h \in I^3$.

We need the following.

\begin{claim*} $h \not\in I$ but $h \in \tilde{I}$. Hence $\tilde{I} \ne I$. 
\end{claim*}

\begin{proof}
If $h \in I$, then $b_1h \in QI$, so that $z_1=b_1h+h'+\tau_3 \in QI,$
whence $f_1=0$ in $S$ (cf. Lemma \ref{fact1} (2)), which is impossible.
Let $1 \leq i \leq d-1$. Then  
$$b_iz_1 - b_1z_i = b_i(b_1h + h' + \tau_3) -b_1z_i =b_1(b_ih -z_i) + b_i(h' + \tau_3) \in Q^2I.$$ Therefore, because $b_i(h' + \tau_3) \in Q^2I$, we get $$b_1(b_ih- z_i) \in (b_1) \cap Q^2I.$$ Notice that 
\begin{eqnarray*}
(b_1) \cap Q^{2}I&=&(b_1) \cap [b_1QI + (b_2, b_3, \cdots, b_d)^{2}I]\\
&=&b_1QI+[(b_1)\cap(b_2,b_3,\cdots,b_d)^{2}I] \\
&=&b_1QI+b_1(b_2,b_3,\cdots,b_d)^{2}\\
&=&b_1QI
\end{eqnarray*}
and we have   
$b_ih- z_i \in QI$, whence $b_ih \in I^2$ for $1 \leq i \leq d-1$. Consequently $b_i^2h \in I^{3}$ for all $1 \leq i \leq d$, so that  $h \in \tilde{I}$, whence $\tilde{I} \ne I$. 
\end{proof}

Because $\ell_A(\tilde{I}/I) \geq 1$, we have
\begin{eqnarray*}
e_1&=& e_0 -\ell_A(A/I)+1\\
   &=& \e_0(\tilde{I})-\ell_A(A/\tilde{I})+[1-\ell_A(\tilde{I}/I)]\\
   & \leq & \e_0(\tilde{I})-\ell_A(A/\tilde{I})\\
   & \leq & \e_1(\tilde{I})\\
   &=& e_1,
\end{eqnarray*}
where $\e_0(\tilde{I}) - \ell_A(A/\tilde{I}) \leq \e_1(\tilde{I})$
is the inequality of Northcott  for the ideal $\tilde{I}$ (cf. Corollary \ref{huneke}).
Hence $\ell_A(\tilde{I}/I)=1$ and $\e_1(\tilde{I})=\e_0(\tilde{I})-\ell_A(A/\tilde{I})$,
so that $$\tilde{I}=I+(h) \ \ \operatorname{and} \ \ \tilde{I}^2=Q \tilde{I}$$ by Corollary \ref{huneke}
(recall that $Q$ is a reduction of $\tilde{I}$ also). We then have, thanks to Proposition \ref{2.6}, that $I^3 = QI^2$, which is a required  contradiction. 
This completes the proof of Theorem  \ref{MainTheorem}  and that of Theorem \ref{I3=QI2} as well. 
\end{proof}

\section{Consequences}

In this section let us review three results of \cite{GNO, S3, V} in order to see how our Theorem \ref{MainTheorem} works to prove or improve them. Let us begin with Theorem \ref{Sally}.

\begin{proof}[Proof of Theorem \ref{Sally}]
Notice that condition (1) (resp. (2)) in Theorem \ref{Sally} is equivalent to condition (3) (resp. (1)) in Theorem \ref{MainTheorem} with $c = 1$. The implication $(1) \Rightarrow (3)$ and the last assertions of Theorem \ref{Sally}  are contained in Theorem \ref{MainTheorem}. Suppose condition (3) of Theorem \ref{Sally} is satisfied. Then $S = TS_1$ since $I^3=QI^2$, whence $\m S = (0)$ and $\mu_B(S) = 1$ because $\ell_A(S_1)=1$ (recall that $S_1 = I^2/QI$ and $\ell_A(I^2/QI)=1$). Thus condition (2) in Theorem \ref{MainTheorem} is satisfied. 
\end{proof}

The following result is the main result of \cite{GNO}, which is exactly the case $c = 2$ of Theorem \ref{MainTheorem}. We would like to refer the reader to \cite{GNO} for the proof, which can be substantially simplified by Theorem \ref{MainTheorem}.

\begin{cor}[\cite{GNO}, Theorem 1.2]\label{gno}
Suppose that $d \geq 2$. Then the following four conditions are equivalent to each other.
\begin{itemize}
\item[$(1)$] $\m S = (0)$, $\operatorname{rank}_BS =1$, and $\mu_B(S) = 2$. 
\item[$(2)$] There exists an exact sequence 
$$0 \to B(-2) \to B(-1) \oplus B(-1) \to S \to 0$$of graded $T$-modules.
\item[$(3)$] $e_1 = e_0 - \ell_A(A/I) + 1$, $e_2 =0$, and $\operatorname{depth} G \geq d-2$.
\item[$(4)$] $I^3 =QI^2$, $\ell_A(I^2/QI) = 2$, $\m I^2 \subseteq QI$, and $\ell_A(I^3/Q^2I) < 2d$.
\end{itemize}
When this is the case, the following assertions hold true
\begin{itemize}
\item[(i)] $\operatorname{depth} G = d-2$.
\item[(ii)] $e_3 = -1$, if $d \geq 3$.
\item[(iii)] $e_i = 0$ for all $4 \leq i \leq d$.
\item[(iv)] $\ell_A(I^3/Q^2I) = 2d - 1$.
\end{itemize}
\end{cor}

Later we need the following result in Section 5, which is due to \cite{GNO} and is exactly the case $c = d$ of Theorem \ref{MainTheorem}. Here we have deleted from the original statement the superfluous condition that $I^3=QI^2$ in conditions (2) and (3) (cf. Proposition \ref{2.6} also). We refer the reader to \cite{GNO} for the proof.

\begin{cor}[\cite{GNO}, Corollary 2.6]\label{B_+}
Suppose that $d \geq 2$. Then the following three conditions are equivalent to each other.
\begin{enumerate}
\item[$(1)$] $S \cong B_+$ as graded $T$-modules.
\item[$(2)$] $e_1 = e_0 - \ell_A(A/I) + 1$ and $e_i = 0$ for all $2 \leq i \leq d$. 
\item[$(3)$] $\ell_A(\tilde{I}/I) =1$ and $\tilde{I}^2 ~= Q\tilde{I}$. 
\end{enumerate}
When this is the case, the graded rings $G$, $R$, and $R'$ are all Buchsbaum rings with Buchsbaum invariant
$${\Bbb I}(G)={\Bbb I}(R)={\Bbb I}(R')=d.$$ 
\end{cor}

We have learned the following example from Rossi.

\begin{ex}
Let $A$ be a $3$-dimensional regular local ring and let $x, y, z$ be  a regular system of parameters. We put  $$I = (x^2 - y^2, x^2 - z^2, xy, yz, zx) \ \ \operatorname{and} \ \ Q = (x^2 - y^2, x^2 - z^2, yz).$$ Then $\tilde{I} = \fkm^2 = I + (z^2)$ and $\ell_A(\tilde{I}/I)= 1$. Since $\fkm^4 = Q\fkm^2$, the ideal $I$ satisfies condition (3) in Corollary \ref{B_+}, so that $\e_1(I) = \e_0(I) - \ell_A(A/I) + 1$ and $\e_2(I) = \e_3(I) =0$. The graded rings $G={\rmG} (I)$, $R = {\rmR}(I)$, and $R'= {\rmR}'(I)$ are all Buchsbaum rings with ${\Bbb I}(G)={\Bbb I}(R)={\Bbb I}(R')=3.$
\end{ex}

\section{An example}
In this section we construct one example which satisfies condition (3) in Theorem \ref{MainTheorem}. Our goal is the following.

\begin{thm}
Let $0 < c \leq d$ be integers. Then there exists an $\fkm$-primary ideal $I$  in a Cohen-Macaulay local ring $(A, \fkm)$ such that  
 $$d = \operatorname{dim}A, \ \ \e_1(I) = \e_0(I) - \ell_A(A/I) + 1, \ \ \operatorname{and} \ \ c = \ell_A(I^2/QI)$$
for some reduction $Q = (a_1, a_2, \cdots, a_d) $ of $I$.
\end{thm}

To construct necessary examples we may assume that $c = d$. In fact, suppose that $0 < c < d$ and assume that we have already chosen an $\fkm_0$-primary ideal $I_0$  in a certain Cohen-Macaulay local ring $(A_0, \fkm_0)$ such that $c = \operatorname{dim}A_0, \ \ \e_1(I_0) = \e_0(I_0) - \ell_{A_0}(A_0/I_0) + 1$, and $c = \ell_{A_0}(I_0^2/Q_0I_0)$ with $Q_0 = (a_1, a_2, \cdots, a_c)A_0$ a reduction of $I_0$. 
Let $n = d - c$ and let $A = A_0[[X_1, X_2, \cdots, X_{n}]]$ be the formal power series ring. We put $I = I_0A + (X_1, X_2, \cdots, X_{n})A$ and $Q = Q_0A + (X_1, X_2, \cdots, X_{n})A$. Then $A$ is a Cohen-Macaulay local ring with $\operatorname{dim} A = \operatorname{dim} A_0 + n = d$ and the maximal ideal $\fkm = \fkm_0A + (X_1, X_2, \cdots, X_{n})A$. The ideal $Q$ is a reduction of $I$ and because $X_1, X_2, \cdots, X_n$ forms a super regular sequence in $A$ with respect to $I$ (recall that ${\rmG} (I) = {\rmG} (I_0)[Y_1, Y_2, \cdots, Y_n]$ is the polynomial ring, where $Y_i$'s are the initial forms of $X_i$'s), we have $\e_i(I) = \e_i(I_0)$~($i = 0, 1$) and $I^2/QI \cong I_0^2/Q_0I_0$, whence $\e_1(I) = \e_0(I) - \ell_A(A/I) + 1$ and $\ell_A(I^2/QI ) = c$. This observation allows us to concentrate our attention on the case where $c = d$.

Let $m,\,d>0$ be integers.
Let $$U=k[\{X_j \}_{1 \leq j \leq m},Y,\{V_i\}_{1 \leq i \leq d}, \{Z_i\}_{1 \leq i \leq d}]$$
be the polynomial ring with $m+2d + 1$ indeterminates over an infinite field $k$ and let 
\begin{eqnarray*}
\a&=&[(X_j \mid 1 \leq j \leq m)+(Y)]{\cdot}[(X_j \mid 1 \leq j \leq m)+(Y)+(V_i \mid 1 \leq i \leq d)]\\
    &&+(V_iV_j \mid  1 \leq i, j \leq d, \,  i \ne j)+(V_i^2-Z_iY  \mid  1 \leq i \leq d).
\end{eqnarray*}
We put $C=U/\a$ and denote the images of $X_j$, $Y$, $V_i$, and $Z_i$ in $C$
by $x_j$, $y$, $v_i$, and $a_i$, respectively.
Then ${\dim}~C=d$, since $\sqrt{\a}= (X_j \mid 1 \leq j \leq m) + (Y) + (V_i \mid 1 \leq i \leq d)$.
Let 
$M = C_+ := (x_j \mid 1 \leq j \leq m) + (y) + (v_i \mid 1 \leq i \leq d) + (a_i \mid 1 \leq i \leq d)$
be the graded maximal ideal in $C$. 
Let $\Lambda$ be a subset of $\{1,2,\cdots,m \}$. We put 
$$J=(a_i \mid  1 \leq i \leq d)+(x_{\alpha} \mid \alpha \in \Lambda)+(v_i \mid 1 \leq i \leq d)~ ~\ \textup{and}~~\ \q=(a_i \mid  1 \leq i \leq d).$$
Then $M^2={\q} M$, $J^2 = {\q}J + {\q}y$, and $J^3={\q}J^2$, 
whence $\q$ is a reduction of both $M$ and $J$, and 
$a_1,a_2,\cdots,a_d$ is a homogeneous system of parameters for the graded ring $C$.

Let $A=C_M$, $I=JA$, and $Q={\q}A$.
We are now interested in the Hilbert coefficients $e_i's$ of the ideal $I$ 
as well as the structure of the associated graded ring and the Sally module of $I$. Let us maintain the same notation as in the previous sections. We then have the following, which shows that the ideal $I$ is a required example.

\begin{thm}\label{ex1}
The following assertions hold true.
\begin{itemize}
\item[$(1)$] $A$ is a Cohen-Macaulay local ring with $\dim A = d$.
\item[$(2)$] $S \cong B_+$ as graded $T$-modules, whence $\ell_A(I^2/QI) = d$.
\item[$(3)$] ${\e}_0(I)=m+d+2$ and ${e}_1=\sharp \Lambda +d+1$.
\item[$(4)$] ${\e}_i(I)=0$ for all $2 \leq i \leq d$.
\item[$(5)$] $G$ is a Buchsbaum ring with  $\operatorname{depth}G=0$ and ${\Bbb I}(G)=d$.
\end{itemize}
\end{thm}

We divide the proof of Theorem \ref{ex1} into a few steps. 
Let us  begin with the following. 

\begin{prop}\label{lemmaF}
Let $\fkp = (X_j \mid  1 \leq j \leq m)+(Y)+(V_i \mid  1 \leq i \leq d)$ in $U$. 
Then ${\ell}_{C_{\p}}(C_{\p})=m+d+2$.
\end{prop}

\begin{proof}
Let $\widetilde{U}=U[\{\frac{1}{Z_i}\}_{1 \leq i \leq d}]$ and put
$\widetilde{k}=k[\{Z_i\}_{1 \leq i \leq d}, \{\frac{1}{Z_i}\}_{1 \leq i \leq d}]$ in $\tilde{U}$.
Let  
$X_j'=\frac{X_j}{Z_1}$~($1 \leq j \leq m$),
$V_i'=\frac{V_i}{Z_1}$~($1 \leq i \leq d$), 
and $Y'=\frac{Y}{Z_1}$.
Then $\{X_j'\}_{1 \leq j \leq m}$, $Y'$, and $\{V_i'\}_{1 \leq i \leq d}$
are algebraically independent over $\widetilde{k}$, 
$$\widetilde{U}=\widetilde{k}[\{X_j'\}_{1 \leq j \leq m}, Y',\{V_i'\}_{1 \leq i \leq d}],\ \ \textup{and}$$
\begin{eqnarray*}
\a\widetilde{U}&=&[(X_j' \mid 1 \leq j \leq m) +(Y')]{\cdot}[(X_j' \mid 1 \leq j \leq m) + (Y') + (V_i' \mid 1 \leq i \leq d)]\\&+&(V_i'V_j' \mid 1 \leq i, j \leq d, \, i \ne j) + (\frac{A_1}{A_i} {V_i'}^2-Y' \mid 1 \leq i \leq d).
\end{eqnarray*}
Let $W=
\widetilde{k}[\{X_j'\}_{1 \leq j \leq m}, \{V_i'\}_{1 \leq i \leq d}]$ in $\widetilde{U}$ and
\begin{eqnarray*}
{\b}&=&[(X_j' \mid 1 \leq j \leq m) + ({V_1'}^2)]{\cdot}[(X_j' \mid 1 \leq j \leq m) + (V_i' \mid 1 \leq i \leq d)]\\
&+&(V_i'V_j' \mid 1 \leq i, j \leq d, \, i \ne j) + (\frac{A_1}{A_i} {V_i'}^2-{V_1'}^2 \mid 2 \leq i \leq d)
\end{eqnarray*}
in $W$. Then, substituting $Y'$ with ${V_1'}^2$ in $\tilde{U}$, we get the isomorphism
$$ \widetilde{U}/\a \widetilde{U} \cong \overline{U}:=W/\b$$
of $\tilde{k}$-algebras, under which the prime ideal ${\p}\widetilde{U}/\a\widetilde{U}$ corresponds to the prime ideal
$P/\b$ of $\overline{U}$, where $P =W_+ := (X_j' \mid 1 \leq j \leq m) + (V_i' \mid 1 \leq i \leq d)$. Then, since $\b+({V_1'}^2)=P^2$ and $$\ell_{W_P}([\b+({V_1'}^2)]W_P/\b W_P)=1,$$
we get 
\begin{eqnarray*}
\ell_{\overline{U}_P}(\overline{U}_P)&=&\ell_{W_P}(W_P/P^2W_P) +\ell_{W_P}([\b+({V_1'}^2)]W_P/\b W_P)\\
&=&(m+d+1) + 1 \\
&=& m + d + 2.
\end{eqnarray*}
Thus $\ell_{C_{\p}}(C_{\p})=\ell_{\overline{U}_{P}}(\overline{U}_{P})=m+d+2$.
\end{proof}

We have by the associative formula of multiplicity that  $${\e}_0({\q})=\ell_{C_{\p}}(C_{\p}) {\cdot} {\e}_0^{C/\fkp C}([\q + \fkp C]/\fkp C)=m+d+2,$$
because $\fkp=\sqrt{\a}$ and $C/\fkp C = U/\fkp =  k[Z_i \mid 1 \leq i \leq d]$.
On the other hand, we have $\ell_C(C/{\q})=m+d+2,$ since $$C/{\q} = k[\{X_j\}_{1 \leq j \leq m}, Y, \{V_i\}_{1 \leq i \leq d}]/\c^2$$
where $$\c = (X_j \mid 1 \leq j \leq m) +(Y) + (V_i \mid 1 \leq i \leq d).$$
Thus  ${\e}_0({\q})=\ell_C(C/{\q})$, so that $C$ is a Cohen-Macaulay ring and $\e_0({\q}) = m+d+2$.

\begin{prop}\label{lemmaex}
$\ell_C(\tilde{J}/J)=1$ and $\tilde{J}^2={\q}\tilde{J}$.
\end{prop}

\begin{proof}
Let $K=J+(y)$. Then  $\ell_C(K/J) = 1$ and $K^2 = \q K = J^2$. Hence $\tilde{K} = K$ because $K^2 =\q K$, while we have $\tilde{K} = \tilde{J}$ because $K^2 = J^2$.  Thus the assertions follow. 
\end{proof}

We are now in a position to finish the proof of Theorem \ref{ex1}.

\begin{proof}[Proof of Theorem \ref{ex1}.]
Since $\tilde{I} = \tilde{J}A$, by Proposition  \ref{lemmaex} we get $\ell_A(\tilde{I}/I)=1$ and $\tilde{I}^2 = Q\tilde{I}$. Hence by Corollary \ref{B_+} $S \cong B_+$ as graded $T$-modules, 
so that $$\e_1(I) = \e_0(I) - \ell_A(A/I) + 1$$ and $\e_i(I) = 0$ for all $2 \leq i \leq d$. We have  $\e_1(I) = \sharp{\Lambda}+d + 1$, because $\ell_A(A/I) = m - \sharp{\Lambda} +2$ and $\e_0(I) = \e_0(\q) = m +d+ 2$.
The ring $G = {\rmG} (I)$ is a Buchsbaum ring with $\operatorname{depth}G=0$ and ${\Bbb I}(G)=d$ by Corollary \ref{B_+}, which completes the proof of Theorem \ref{ex1}.
\end{proof}


\addcontentsline{toc}{section}{references}

\begin{thebibliography}{GNO}



\bibitem[CPP]{CPP} A. Corso, C. Polini and M. Vaz Pinoto, \textit{Sally modules and associated graded rings}, Comm. Algebra, $\bf26$, 1998, 2689--2708.


\bibitem[GNO]{GNO} S. Goto, K. Nishida, and K. Ozeki, \textit{Sally modules of rank one}, Michigan Math. J. (to appear).


\bibitem[H]{H} C. Huneke, \textit{Hilbert functions and symbolic powers}, Michigan Math. J., Vol $\bf34$, 1987, 293--318.




\bibitem[N]{N} D. G. Northcott, \textit{A note on the coefficients of the abstract Hilbert function}, J. London Math. Soc., Vol $\bf35$, 1960, 209--214.


\bibitem[RR]{RR} L. J. Ratliff and D. Rush, \textit{Two notes on reductions of ideals}, Indiana Univ. Math. J. $\bf27$, 1978, 929-934.

\bibitem[R]{R} M. E. Rossi, \textit{A bound on the reduction number of a primary ideal}, Proc. Amer. Math. Soc. 128, 2000, no. $\bf5$, 1325--1332.


\bibitem[S1]{S1} J. D. Sally, \textit{Cohen-Macaulay local rings of maximal embedding dimension}, J. Algebra, $\bf56$, 1979, 168--183.

\bibitem[S2]{S2} J. D. Sally, \textit{Tangent cones at Gorenstein singularities}, Composito Math. $\bf40$, 1980, 167--175.

\bibitem[S3]{S3} J. D. Sally, \textit{Hilbert coefficients and reduction number 2}, J. Alg. Geo. and Sing. $\bf1$, 1992, 325--333.



\bibitem[V]{V} W. V. Vasconcelos, \textit{Hilbert Functions, Analytic Spread, and Koszul Homology}, Contemporary Mathematics, Vol $\bf159$, 1994, 410--422.


\bibitem[W]{W} H.-J. Wang, 
{\it Links of symbolic powers of prime ideals,}
Math. Z., {\bf 256} (2007), 749--756.




\bibitem[P]{P} M. Vaz Pinoto, \textit{Hilbert functions and Sally modules}, J. Algebra, $\bf192$, 1977, 504--523.
\end{thebibliography}
\end{document}